\documentclass[11pt,a4paper]{article}
\usepackage{amsmath,amssymb,amsfonts,amsthm}

\newcommand{\E}{{\bf{E}}}
\newcommand{\PP}{{\bf{P}}}
\newcommand{\Var}{{\bf{Var}}}

\newtheorem{thm}{Theorem}

\newtheorem{lem}{Lemma}
\newtheorem{prop}{Proposition}
\theoremstyle{definition}

\begin{document}

\bibliographystyle{plain}
\parindent=0pt
\centerline{\LARGE \bfseries  The largest component}
\centerline{\LARGE \bfseries in an inhomogeneous random intersection graph}

\centerline{\LARGE \bfseries with clustering}
{\large  \qquad\qquad\qquad \qquad \qquad}\\

\par\vskip 2.5em

\centerline{ M. Bloznelis}

\bigskip

 Faculty of Mathematics and Informatics, Vilnius University,
LT-03225 Vilnius, Lithuania

E-mail \ {\it mindaugas.bloznelis@mif.vu.lt}


\par\vskip 1.5em

\centerline{February 20, 2010}

\par\vskip 1.5em
\begin{abstract}Given integers $n$, $m=\lfloor\beta n \rfloor$ and a probability measure $Q$
on $\{0, 1,\dots, m\}$, consider the random intersection graph
on the vertex set $[n]=\{1,2,\dots, n\}$, where $i,j\in [n]$ are declared adjacent whenever
$S(i)\cap S(j)\not=\emptyset$.
Here $S(1),\dots, S(n)$
denote iid random subsets
of $[m]$ with the distribution $\PP(S(i)=A)={{m}\choose{|A|}}^{-1}Q(|A|)$,  $A\subset [m]$.
For sparse random intersection graphs  we establish a first order asymptotic as $n\to \infty$
for the order of the largest connected component $N_1=n(1-Q(0))\rho+o_P(n)$.
Here $\rho$ is an average of  nonextinction probabilities of
a related multi-type Poisson branching process.
\par\end{abstract}

\section{Introduction}

Let $Q$ be a probability measure on $\{0,1,\dots, m\}$, and let  $S_1,\dots, S_n$ be random subsets of a set
$W=\{w_1,\dots, w_m\}$ drawn independently from the probability distribution
$\PP(S_i =A)={{m}\choose{|A|}}^{-1}Q(|A|)$,  $A\subset W$, for $i=1,\dots, n$.
A random intersection graph $G(n,m,Q)$ with a vertex set $V=\{v_1,\dots, v_n\}$ is defined as follows.
Every vertex $v_i$ is prescribed the set $S(v_i)=S_i$ and two vertices $v_i$ and $v_j$ are declared adjacent
(denoted $v_i\sim v_j$)
whenever $S(v_i)\cap S(v_j)\not=\emptyset$.
The elements of $W$ are sometimes called attributes, and $S(v_i)$ is called the set of attributes of $v_i$.

Random intersection graphs $G(n,m,Q)$ with the binomial distribution $Q\sim Bi(m,p)$ were
introduced in Singer-Cohen \cite{singerphd}
and Karo\'nski et al. \cite{karonski1999}, see also  \cite{fill2000} and \cite{stark2004}.
The emergence of a giant connected component
in a sparse binomial random intersection graph was studied by Behrish \cite{behrisch2005},
for $m=\lfloor n^{\alpha} \rfloor$, $\alpha\not=1$, and by  Lager{\aa}s and Lindholm \cite{Lageras}, for
  $m=\lfloor \beta n\rfloor$, where $\beta>0$ is a constant.
   They have shown, in particular, that,  for $\alpha\ge 1$, the largest connected component  collects a fraction of all
    vertices whenever the
   average vertex degree, say $d$,
   is larger than $1+\varepsilon$. For $d<1-\varepsilon$ the order of the
   largest connected component is $O(\log n)$.

The graph $G(n,m,Q)$ defined by an arbitrary probability measure $Q$ (we call such graphs inhomogeneous) was first
considered
in Godehardt and Jaworski \cite{godehardt2003}, see also  \cite {JKS}.  Deijfen and Kets \cite{Deijfen}, and
 Bloznelis
\cite{Bloznelis2007/2} showed (in increasing generality) that the typical vertex degree of $G(n,m,Q)$ has the  power law for a heavy tailed distribution
$Q$.  Another result by Deijfen and Kets \cite{Deijfen} says that, for $m\approx \beta n$, graphs $G(n,m,Q)$ posses
the clustering property.

The emergence of a giant connected component
in a sparse  inhomogeneous intersection graph with $n=o(m)$ (graph without clustering) was studied in
\cite{Bloznelis2008+}.
The present paper addresses
inhomogeneous intersection graphs with
clustering, i.e., the case where $m\approx \beta n$.

\section{Results}

Given $\beta>0$, let  $\{G(n,m_n,Q_n)\}$ be a sequence of random intersection graphs such that
\begin{equation}\label{m-n}
\lim_nm_nn^{-1}=\beta.
\end{equation}
 We shall assume that
the sequence of probability distributions $\{Q_n\}$ converges to some probability distribution
$Q$ defined on $\{0,1,2,\dots\}$,
\begin{equation}\label{Teorema-1}
\lim_nQ_n(t)=Q(t),
\qquad
\forall \ t=0,1,\dots,
\end{equation}
and, in addition, the sequence of the first moments converges,
\begin{equation}\label{Teorema-2}
\lim_{n}\sum_{t\ge 1}tQ_n(t)=\sum_{t\ge 1}tQ(t)<\infty.
\end{equation}

{\bf 2.1. Degree distribution.} Let $V_n=\{v_1,\dots, v_n\}$
 denote the vertex set
 of  $G_n=G(n,m_n,Q_n)$ and let $d_n(v_i)$ denote the degree
 of vertex $v_i$.
 Note that, by symmetry, the random variables $d_n(v_1),\dots, d_n(v_n)$ have the same probability
 distribution, denoted $D_n$. In the following proposition we recall a known fact about the asymptotic distribution of $D_n$.
 \begin{prop}
 Assume that (\ref{m-n}), (\ref{Teorema-1}) and (\ref{Teorema-2}) hold. Then we have as $n\to\infty$
 \begin{equation}\label{degree}
 \PP(D_n=k)\to \sum_{t\ge 0}\frac{(at)^k}{k!}e^{-at}Q(t),
 \qquad
 k=0,1,\dots .
 \end{equation}
 Here $a=\beta^{-1}\sum_{t\ge 0}tQ(t)$.
 \end{prop}
 Roughly speaking, the limiting distribution of $D_n$ is the Poisson distribution ${\cal P}(\lambda)$ with random
 parameter $\lambda=aX$,
 where $X$ is a random variable with the distribution $Q$. In particular, for a heavy tailed distribution $Q$ we obtain
 the heavy tailed asymptotic distribution
 for $D_n$.
  For  $Q\sim Bi(m,p)$,  (\ref{degree}) is shown in
 \cite{stark2004}.  For arbitrary $Q$, (\ref{degree}) is shown (in increasing generality)  in \cite{Deijfen} and \cite{Bloznelis2009a}.

 \smallskip

 {\bf 2.2. The largest component.}
Let $N_1(G)$ denote the
order of the largest connected component of a graph $G$ (
i.e., $N_1(G)$ is the number of vertices of a connected component
which has the largest number of vertices). We are interested in a first order asymptotic
of $N_1(G(n,m_n,Q_n))$ as $n\to\infty$.

The most commonly used approach to the parameter $N_1(G)$  of a random graph $G$ is based
on tree counting,
see \cite{erdos1960}, \cite{Chung}. For inhomogeneous random graphs  it is convenient to count trees with a help
of branching processes, see \cite{Riordan}. Here large trees correspond to surviving branching processes and the order
of the largest connected component is described by means of the survival probabilities of a related branching process.

In the present paper we use the approach developed in  \cite{Riordan}. Before formulating our main  result Theorem
\ref{theorem} we will introduce some notation.
Let ${\cal X}={\cal X}_{Q,\beta}$ denote the multi-type
Galton-Watson
branching process, where particles are of
types $t\in {\mathbb T}=\{1,2,\dots\}$ and where the number of
children of type $t$ of a particle of type $s$ has the Poisson distribution with mean $(s-1)tq_t\beta^{-1}$.
Here we write $q_t=Q(t)$,
$t\in {\mathbb T}$.
Let ${\cal X}(t)$ denote the process ${\cal X}$ starting at a particle of type $t$, and $|{\cal X}(t)|$ denote
the total progeny of ${\cal X}(t)$. Let $\rho_{Q,\beta}(t)=\PP(|{\cal X}(t)|=\infty)$ denote the survival probability of the
process ${\cal X}(t)$. Write
$\rho^{(k)}_{Q,\beta}(t)=\PP(|{\cal X}(t)|\ge k)$,
\begin{displaymath}
{\tilde \rho}_{Q,\beta}=\sum_{t\in {\mathbb T}}\rho_{Q,\beta}(t+1)q_t,
\qquad
{\tilde \rho}^{(k)}(Q)=\sum_{t\in {\mathbb T}}\rho^{(k)}_{Q,\beta}(t+1)q_t.
\end{displaymath}
Note that for every $t\in {\mathbb T}$ we have $\rho^{(k)}_{Q,\beta}(t)\downarrow\rho_{Q,\beta}(t)$
as $k\uparrow\infty$ (by the continuity property of  probabilities).
Hence, ${\tilde \rho}^{(k)}(Q)\downarrow{\tilde \rho}(Q)$
as $k\uparrow\infty$.

\begin{thm}\label{theorem}
Let $\beta>0$. Let  $\{m_n\}$ be a sequence of integers satisfying (\ref{m-n}).
  Let
$Q, Q_1, Q_2,\dots$
be probability measures
defined on $\{0,1,2\dots\}$ such that $\sum_{t=0}^{m_n}Q_n(t)=1$, for $n=1,2,\dots$.
Assume that (\ref{Teorema-1}) and (\ref{Teorema-2}) holds.
 Then we have as $n\to\infty$
\begin{equation}\label{Teorema-3}
N_1\bigl(G(n,m_n,Q_n)\bigr)=n\bigl(1-Q(0)\bigr)\rho +o_P(n).
\end{equation}
Here $\rho={\tilde \rho}_{Q^*,\beta^*}$, for $Q(0)<1$, and $\rho=0$, otherwise.
 $Q^*$ denotes the probability measure on $\{1,2,\dots\}$
 defined by $Q^*(t)=\bigl(1-Q(0)\bigr)^{-1}Q(t)$, $t\ge 1$, and $\beta^*=\beta\bigl(1-Q(0)\bigr)^{-1}$.
\end{thm}

Notation $o_P(n)$. We write $\eta_n=o_P(1)$ for a sequence of random variables
$\{\eta_n\}$ that converges to $0$ in probability.
We write  $\eta_n=o_P(n)$ in the case where $\eta_nn^{-1}=o_P(1)$.

{\it Remark 1}. The correspondence $\rho>0\Leftrightarrow \E D_n>c>1$ established for binomial random intersection graphs
in \cite{behrisch2005}, \cite{Lageras} can not be extended to  general inhomogeneous graphs $G(n,m_n,Q_n)$. To see this,
consider the graph obtained from a binomial random intersection graph by replacing $S(v_i)$ by $\emptyset$ for
a randomly chosen fraction of vertices.  This way we can make the expected degree arbitrarily small, and still have
the giant connected component spanned by a fraction of unchanged vertices.

{\it Remark 2}. The kernel $(s,t)\to (s-1)t\beta^{-1}$ of the Poisson branching process which determines the fraction
$\rho$ in the case $m_n\approx \beta n$ differs from  the kernel $(s,t)\to st$ which appears in the case $n=o(m_n)$, see
\cite{Bloznelis2008+}.

\section{Proof}
The section is organized as follows. Firstly we collect some notation and formulate auxiliary results. We then prove Theorem
\ref{theorem}. The proofs of auxiliary results are given in the end of the section.

\smallskip

Let $W'$ be a finite set of size $|W'|=k$. Let $B,H$ be subsets of $W'$ of sizes $|B|=b$ and $|H|=h$ such that
$B\cap H=\emptyset$. Let $A$ be a random subset of $W'$ uniformly distributed in the class of subsets of $W'$ of size $a$.
Introduce the probabilities
\begin{eqnarray}
\nonumber
&&
p(a,b,k)=\PP(A\cap B\not=\emptyset),
\\
&&
\nonumber
p_1(a,b,k)=\PP(|A\cap B|=1),
\qquad
p_2(a,b,k)=\PP(|A\cap B|\ge 2),
\\
\nonumber
&&
p(a,b,h,k)=\PP\bigl(|A\cap B|=1,\, A\cap H=\emptyset\bigr),
\\
\nonumber
&&
p_1(a,b,h,k)=\PP\bigl(|A\cap B|=1,\, A\cap H\not=\emptyset\bigr).
\end{eqnarray}

\begin{lem}\label{inequalities-1} Let $k\ge 4$.  Denote $\varkappa=ab/k$ and $\varkappa'=ab/(k-a)$.
For $a+b\le k$ we have
\begin{eqnarray}
\label{pab}
\varkappa(1-\varkappa')
\le
&
p_1(a,b,k)
&
\le p(a,b,k)\le \varkappa,
\\
\label{p2ab}
&
p_2(a,b,k)
&
\le 2^{-1}\varkappa^2.
\end{eqnarray}
Denote $\varkappa''=(a-1)h/(k-b)$. For $a+b+h\le k$ we have
\begin{eqnarray}
\label{pabh}
\varkappa(1-\varkappa'-\varkappa'')
\le
&
p(a,b,h,k)
&
\le \varkappa.
\\
\label{pabh+1}
&
p_1(a,b,h,k)
&
\le \varkappa_h\varkappa.
\end{eqnarray}
\end{lem}

Given integers $n, m$ and a vector ${\overline s}=(s_1,\dots, s_n)$ with   coordinates from the set
$\{0,1,\dots, m\}$,
let $S(v_1),\dots, S(v_n)$ be independent random subsets of $W_m=\{w_1,\dots, w_m\}$ such that,
for every $1\le i\le n$, the subset $S(v_i)$ is uniformly distributed in the class of all
subsets of $W_m$ of size
$s_i$. Let
$G_{\overline s}(n,m)$ denote the random intersection graph on the vertex set $V_n=\{v_1,\dots, v_n\}$
defined by the random sets $S(v_1),\dots, S(v_n)$. That is, we have $v_i\sim v_j$ whenever
$S(v_i)\cap S(v_j)\not=\emptyset$.

\begin{lem}\label{Discrete}
Let $M>0$ be an integer and let $Q$ be a probability measure defined on
\linebreak
$[M]=\{1,\dots, M\}$.
Let  $\{m_n\}$ be a sequence of integers, and
$\{ {\overline s}_n=(s_{n1},\dots, s_{nn})\}$
be a sequence of vectors with integer coordinates
$s_{ni}\in [M]$, $1\le i\le n$.
 Let
$n_t$ denote the  number of coordinates of ${\overline s}_n$  attaining the value $t$.
Assume that, for some integer $n'$ and a sequence $\{\varepsilon_n\}$ $\subset (0,1)$
converging to zero,
we have, for every $n>n'$,
\begin{eqnarray}\label{2010--A1}
&&
\qquad
 \max_{1\le t\le M}|(n_t/n)-Q(t)|\le \varepsilon_n,
\\
\label{2010--A2}
&&
\qquad
 | m_n(\beta n)^{-1}-1|\le \varepsilon_n.
\end{eqnarray}
Then there exists
a sequence $\{\varepsilon^*_n\}_{n\ge 1}$ converging to zero   such that, for $n>n'$, we have
\begin{equation}\label{Discrete-1}
\PP\bigl(\bigl|N_1(G_{{\overline s}_n}(n,m_n))-n{\tilde \rho}_{Q,\beta}\bigr|>\varepsilon^*_nn\bigr)<\varepsilon^*_n.
\end{equation}
\end{lem}

Several technical steps of the proof of Lemma \ref{Discrete} are collected in the separate Lemma \ref{coupling}.

\begin{lem}\label{coupling} Assume that conditions of lemma \ref{Discrete} are satisfied.
For any function $\omega(\cdot)$ satisfying $\omega(n)\to+\infty$ as $n\to\infty$
bounds (\ref{2010-A2}),  (\ref{2010-AA1}), and (\ref{2010-AA3}) hold true.
\end{lem}

\begin{proof}[Proof of Theorem \ref{theorem}] Write, for short, $G_n=G(n,m_n, Q_n)$ and  $N_1=N_1(G(n,m_n, Q_n))$.
Given $t=0,1,\dots$, let $n_t$ denote the number of vertices of $G$ with the attribute sets
of size $t$. Write $q_{nt}=Q_n(t)$ and $q_t=Q(t)$, and $q^*_t=Q(t)$.

Note that  vertices with empty attribute sets are isolated in $G$. Hence, the connected components of order at least $2$
of $G$ belong to the subgraph $G_{[\infty]}\subset G$ induced by the vertices with non-empty attribute sets.

In the case where   $q_0=1$, we obtain from  (\ref{Teorema-1}) that the expected number of vertices in $G_{[\infty]}$
$\E(n- n_0)=n(1-q_{n0})=o(n)$. This identity implies $N_1=o_P(n)$. We obtain (\ref{Teorema-3}), for $q_0=1$.

Let us prove (\ref{Teorema-3}) for $q_0<1$.
Let $G_{[M],n}$ denote the subgraph
of $G_n$ induced by the vertices with attribute sets of sizes from the set $[M]$.
In the proof we approximate $N_1(G_n)$ by $N_1(G_{[M],n})$
and use the result
for $N_1(G_{[M],n})$  shown in
 Lemma \ref{Discrete}.

We need some notation related to $G_{[M],n}$. The inequality $q_0<1$ implies that,
for large $M$,
the sum $q_{[M]}:=q_1+\dots+q_M\approx 1-q_0$ is positive.
Given such $M$, let $Q^*_M$ be the probability measure  on
$[M]$, which assigns the mass $q^*_{M t}=q_t/q_{[M]}$ to  $t\in [M]$.
Denote ${\tilde \rho}_{[M]}={\tilde \rho}_{Q^*_M,\beta_M}$, where $\beta_M=\beta/q_{[M]}$.
Clearly, $\beta_M$ converges to $\beta^*$ as $M\to\infty$, and we have
\begin{equation}\label{M}
\forall t\ge 1 \quad \lim_Mq^*_{Mt}=q^*_t
\qquad
{\text{and}}
\qquad
\lim_M\sum_{t\ge 1}tq^*_{Mt}=\sum_{t\ge 1}tq^*_t<\infty.
\end{equation}
 It follows from (\ref{M}) that
\begin{equation}\label{B-M}
\lim_M{\tilde\rho}_{[M]}={\tilde \rho}_{Q^*,\beta^*}.
\end{equation}
For the proof of (\ref{B-M}) we refer to Chapter 6 of \cite{Riordan}.

We are now ready to prove (\ref{Teorema-3}). For this purpose we combine the upper and lower bounds
\begin{equation}\nonumber
N_1\ge n(1-q_0){\tilde \rho}_{Q^*,\beta^*} -o_P(n)
\quad
{\text{and}}
\quad
N_1\le n(1-q_0){\tilde \rho}_{Q^*,\beta^*} +o_P(n).
\end{equation}
We give the proof of the lower bound only. The proof of the upper bound is almost the same as that of a corresponding
bound in \cite{Bloznelis2008+},
see  formula (56) in \cite{Bloznelis2008+}.

In the proof
we show that, for every $\varepsilon\in (0,1)$,
\begin{equation}\label{vasario-17+}
\PP(N_1>n(1-q_0){\tilde \rho}_{Q^*,\beta^*}-2\varepsilon n)=1-o(1)
\quad
{\text{as}}
\quad
n\to\infty.
\end{equation}
Fix $\varepsilon\in (0,1)$. In view of (\ref{B-M}) we can choose $M$ such that
\begin{equation}\label{rrr-16}
{\tilde \rho}_{Q^*,\beta^*}-\varepsilon<{\tilde \rho}_{[M]}< {\tilde \rho}_{Q^*,\beta^*}+\varepsilon.
\end{equation}
We apply Lemma \ref{Discrete} to $G_{[M],n}$ conditionally given the event
\begin{displaymath}
{\cal A}_n=\{\max_{1\le t\le M}|n_t-q_tn|<n\delta_n+n^{2/3}\}.
\end{displaymath}
Here  $\delta_n=\max_{1\le t\le M}|q_{nt}-q_t|$ satisfies  $\delta_n=o(1)$, see (\ref{Teorema-1}). In addition, we have
 \begin{eqnarray}
\nonumber
1-\PP({\cal A}_n)
&&
\le
\PP(\max_{1\le t\le M}|n_t-q_{nt}n|\ge n^{2/3})
\\
\nonumber
&&
\le
\sum_{1\le t\le M}\PP(|n_t-q_{nt}n|\ge n^{2/3})
\\
\nonumber
&&
\le
M\, n^{-1/3}=o(1).
\end{eqnarray}
In the last step we have invoked the  bounds $\PP(|n_t-q_{nt}n|\ge n^{2/3})\le n^{-1/3}$, which follow by
Chebyshev's inequality applied to binomial random variables $n_t$, $t\in [M]$.
Now, combining the bound, which follows from Lemma \ref{Discrete},
\begin{equation}\label{rrr-17}
\PP\bigl(|N_1(G_{[M],n})-n{\tilde \rho}_M|>n\varepsilon\bigr|{\cal A}_n)=o(1)
\end{equation}
with (\ref{rrr-16}) and the bound $\PP({\cal A}_n)=1-o(1)$, we obtain
\begin{displaymath}
\PP\bigl(|N_1(G_{[M],n})-n{\tilde \rho}_{Q^*,\beta^*}|>2n\varepsilon\bigr)=o(1).
\end{displaymath}
Finally, (\ref{vasario-17+}) follows from the obvious inequality $N_1\ge N_1(G_{[M],n})$.

\end{proof}

\begin{proof}[Proof of Lemma \ref{Discrete}]
The proof consists of two steps. Firstly, we show that  components of order at least  $n^{2/3}$  contain
 $n{\tilde \rho}_{Q,\beta}+o_P(n)$ vertices in total. This implies the upper bound for $N_1=N_1(G_{{\overline s}_n}(n,m_n))$
  \begin{equation}\label{upper-bound}
N_1\le n {\tilde \rho}_{Q,\beta}+o_P(n).
\end{equation}
Secondly, we prove that  with probability
tending to one such vertices belong to a common connected component. This implies the lower bound
\begin{equation}\label{lower-bound}
N_1\ge n {\tilde \rho}_{Q,\beta}-o_P(n).
\end{equation}
Clearly, (\ref{upper-bound}), (\ref{lower-bound}) yield (\ref{Discrete-1}).
Before the proof of (\ref{upper-bound}), (\ref{lower-bound}), we
introduce some notation.

{\it Notation.} Denote ${\tilde \rho}={\tilde \rho}_{Q,\beta}$ and write
$q_t=Q(t)$, $t\in [M]$.
  In what follows, we drop the subscript $n$ and write $m=m_n$,
$V=V_n$, $W=W_m$,
$G=G_{{\overline s}_n}(n,m)$.
  We say that a vertex  $v\in V$ is of type $t$  if
the size $s_v=|S(v)|$
of its attribute set $S(v)$ is  $t$.
An edge $u'\sim u''$ of $G$ is called regular if $|S(u')\cap S(u'')|=1$. In this case $u'$ and $u''$ are
called regular neighbours. The edge $u'\sim u''$ is called irregular otherwise.
We say that $v_i$ is smaller than $v_j$
whenever $i<j$.
Given $v\in V$, let $C_v$ denote the connected component of $G$ containing vertex $v$.

In order to count vertices of $C_v$ we explore this component using the  Breath-First Search procedure.

{\it Component exploration.}   Select $v\in V$.
In the beginning all vertices are uncoloured.
Colour $v$ white and add it to the list $L_v$ (now $L_v$ consists of a single white vertex $v$).
%
%
Next we proceed recursively. We choose the oldest white vertex in the list, say $u$,
scan the current set of uncoloured vertices (in increasing order) and look for
 neighbours of $u$.
Each new discovered neighbour  immediately receives white colour and is added to the list.
In particular, neighbours with smaller indices are added to the list before ones with larger indices.
 Once all the uncoloured vertices are scanned colour $u$ black.
 Neighbours of $u$ discovered in this step are called children of $u$.
 We say that $u'\in L_v$ is older than
$u''\in L_v$ if $u'$ has been added to the list before $u''$.
Exploration ends when there are no more white vertices in the list available.

By $L^*_v=\{v=u_1,\, u_2,\, u_3,  \dots\}$ we denote the final state of the list
after the exploration is complete.
Here $i<j$ means that $u_i$ has been discovered before $u_j$.
Clearly, $L^*_v$ is the vertex set of $C_v$.
Denote  $L_v(k)=\{u_i\in L^*_v:\, i\le k\}$. Note that $|L_v(k)|=\min\{k, |L^*_v|\}$.
By $u_{j^*}$ we denote the vertex which has discovered
$u_j$ ($u_j$ is a child of $u_{j^*}$).
Introduce the sets,
\begin{equation}\label{sets-DHS}
D_k=\cup_{1\le j\le k}S(u_j),
\qquad
S'(u_i)=S(u_i)\setminus D_{i-1},
\qquad
k\ge 1,
\quad
i\ge 2,
\end{equation}
and put $D_0=\emptyset$, $S'(u_1)=S(u_1)$.

%
%

{\it Regular}  exploration is performed similarly to
the 'ordinary' exploration, but now only regular neighbours  are added to the
list. We call them regular children. A regular child $u'$ of $u$ is called simple if $S(u')\setminus S(u)$ does not intersect with  $S(e)$
for any vertex $e$ that has  already been included in the list before $u'$. Otherwise the regular child is called
{\it complex}.
{\it Simple} exploration is performed similarly to the {\it regular} exploration, but now  simple
children are added to the list only.

In the case of  regular (respectively simple) exploration we use
the notation $L_v^r$, $L_v^{r*}$, $L_v^r(k)$, $D_k^r$,
${S'}^r(u_i)$
 (respectively $L_v^s$, $L_v^{s*}$, $L_v^s(k)$,
$D_k^s$, ${S'}^s(u_i)$) which is defined
 in much the same way as above. Similarly,
 $i^*$ denotes the number in the list ($L_v^r$ or $L_v^s$ depending on the context) of the vertex
that has discovered
$u_i$ ($u_i$ is a child of $u_{i^*}$).
%
%
%
%
%
 For a member $u_j$ of the list
$L_v^{s*}=\{v=u_1,u_2,\dots\}$ we denote $H (u_j)=(\cup_{j^*<r<j}S(u_r))\setminus D^s_{j^*}$.
Consider the simple exploration at the moment where the current oldest white vertex, say $u_i$ of evolving list
$L_v^s=\{v=u_1, u_2,\dots\}$ starts the search of its simple children. Let
$U_i=\{v_{j_1},\dots,v_{j_r},\dots\, v_{j_k}\}$
denote the current set of uncoloured vertices (the set of potential simple children).
Here $j_1<j_2<\dots<j_k$. Firstly, allow
$u_i$ to discover its simple children among $\{v_{j_1},\dots, v_{j_r-1}\}$.
Define the set $H_i(v_{j_r})=\bigl(\cup_{u\in L}S(u)\bigr)\setminus D^s_i$, where
$L$ denotes the set of current white elements of the list  that are younger than
$u_i$. In particular, $L$ includes the simple children of $u_i$ discovered among $v_{j_1},\dots, v_{j_r-1}$.
 Observe that any
 $u'\in U_i$ becomes a simple child of
$u_i$ whenever
it is a regular neighbour of $u_i$ and $H_i(u')\cap S(u')=\emptyset$.
\begin{equation}
|S(u')\cap S(u_i)|=1
\qquad
\quad
{\text{and}}
\qquad
\quad
S(u')\cap H_i(u')=\emptyset.
\end{equation}
Observe that for
any  member of the list $u_j\in L_v^{s*}$ we have $H(u_j)=H_{j^*}(u_j)$.

 Note that irregular neighbours discovered during {\it regular} exploration  receive white
 colour, but are not added to the list $L_v^r$. Similarly, irregular neighbours and complex children discovered
 during
 {\it simple} exploration receive white colour, but are not added to the list $L_v^s$.
  Note also that $L_v^{s*}$ does not need to be a subset of $L_v^{r*}$.

 Let $\omega(n)$ be an integer function such that $\omega(n)\to+\infty$ and $\omega(n)=o(n)$ as $n\to\infty$.
 A vertex $v\in V$ is called big (respectively, br-vertex and  bs-vertex)  if
 $|L_v^{*}|\ge\omega(n)$ (respectively, $|L^{r*}_v|\ge\omega(n)$ and $|L^{s*}_v|\ge \omega(n)$).
  Let $B$, $B^r$, and $B^s$ denote the collections of big vertices, br-vertices,  and bs-vertices respectively.
 Clearly, we have $B^s,\, B^r\subset B$. Note that in order to decide whether a vertex $v$ is big
  we do not need to explore the component $C_v$ completely. Indeed, we may stop the exploration after the
 number of coloured
 vertices reaches $\omega(n)$.  In what follows we assume that the exploration was stopped after the number
 of coloured vertices had reached $\omega(n)$ (in this case $v\in B$) or ended even earlier because
 the last  white vertex of the list failed to find an uncoloured neighbour (in this case $v\notin B$).

 {\it The upper bound.} Fix $\omega(\cdot)$. We show that
 \begin{equation}\label{2010-1}
 |B|-n{\tilde\rho}=o_P(n).
 \end{equation}
 Note that (\ref{2010-1}) combined with the simple inequality $N_1\le\max\{\omega(n), |B|\}$ implies
(\ref{upper-bound}).
We obtain (\ref{2010-1}) from the bounds
  \begin{eqnarray}\label{2010-A1}
 &&
 |B|-|B^s|=o_P(n),
 \\
 \label{2010-A2}
 &&
 |B^s|-n{\tilde \rho}=o_P(n).
 \end{eqnarray}
(\ref{2010-A2}) is shown in Lemma \ref{coupling}.
(\ref{2010-A1}) follows from the bound
$\E(|B|-|B^s|)=o(n)$. In order to prove this bound we show that
 \begin{eqnarray}\label{2010-AA1}
 &&
 \E|B^s|-n{\tilde \rho}=o(n),
 \\
 \label{2010-AA2}
&&
\E|B| \ \, \le \ n{\tilde \rho}+o(n).
\end{eqnarray}
(\ref{2010-AA1}) is shown in {\bf Lemma } \ref{coupling}.  (\ref{2010-AA2}) follows from the bounds
\begin{eqnarray}\label{2010-AA3}
&&
\E|B^r| \le  n{\tilde \rho}+o(n),
\\
&&
\label{2010-AA4}
\E|B\setminus B^r|=o(n).
\end{eqnarray}
(\ref{2010-AA3}) is shown in {\bf Lemma } \ref{coupling}. In order to show (\ref{2010-AA4})
we write
$\E|B\setminus B^r|=\sum_{v\in V}\PP(v\in B\setminus B^r)$ and invoke the bounds, which hold uniformly in $v\in V$,
 \begin{equation}\label{2010-B}
 \PP(v\in B\setminus B^r)=O(\omega(n)n^{-2}).
 \end{equation}
In the proof of (\ref{2010-B}) we inspect the list $L_v(\omega(n))$ and look for an irregular child.
The probability that given $u_i\in L_v(\omega(n))$ is an
irregular child is $O(n^{-2})$, see (\ref{p2ab}).
Now  (\ref{2010-B}) follows from the fact that $L_v(\omega(n))$ has at most $\omega(n)=o(n)$ elements. The proof of (\ref{2010-A1}) is complete.

\bigskip
{\it The lower bound.}
We start with a simple observation that whp each attribute $w\in W$ is shared by at most $O(\ln n)$ vertices. Denote
 $f(w)=\sum_{v\in V}{\mathbb I}_{\{w\in S(v)\}}$, $w\in W$. We show that the inequality
\begin{equation}\label{f-N}
\max_{w\in W}f(w)\le 2M\ln n
\end{equation}
holds with probability $1-o(1)$.
Since $f(w)$
is
a sum of independent Bernoulli
random variables with success probabilities at most $M/m$, Chernoff's inequality implies
$\PP(f(w)>2M\ln n)\le c_{M,\beta}n^{-2}$. Hence, the complementary event to (\ref{f-N}) has probability
\begin{displaymath}
\PP(\max_{w\in W}f(w)> 2M\ln n)\le \sum_{w\in W}\PP(f(w)>2M\ln n)=o(1).
\end{displaymath}

Let us prove (\ref{lower-bound}).
Fix $\varepsilon\in (0,1)$. For each $t\in [M]$ choose
$\lceil n_t\varepsilon\rceil$
vertices of type $t$
and colour them red. Let $G'$ denote the subgraph of $G$ induced by uncoloured vertices,
 and let $C_1,C_2,\dots$ denote the (vertex sets of) connected components of $G'$ of order at least
$n^{2/3}$. Observe, that the number, say $k$, of such components is at most $(1-\varepsilon)n^{1/3}$.
We apply
(\ref{2010-1}) to the intersection graph $G'$ and function $\omega(n)=\lceil n^{2/3}\rceil$ and obtain
$|\cup_{i\ge 1}C_i|=(1-\varepsilon) n{\tilde\rho}_{Q,\beta'}+o_P(n)$, where $\beta'=\beta(1-\varepsilon)^{-1}$.
We show below that
with a high probability all vertices of  $\cup_{i\ge 1}C_i$ belong to a single connected component of the graph $G$.
Hence,
$N_1\ge (1-\varepsilon){\tilde\rho}_{Q,\beta'}+o_P(n)$.
 Letting $\varepsilon\to 0$ we then immediately obtain lower bound (\ref{lower-bound}).

We assume that $G$ is obtained in two steps. Firstly, the uncoloured vertices generate  $G'$, and, secondly, the red vertices
add the remaining part of $G$. Let us consider the second step where the red vertices add their contribution.
 Write ${\mathbb I}_{ij}=1$  if $C_i$ and $C_j$ are not connected by a path in $G$, and
${\mathbb I}_{ij}=0$ otherwise. Let
 $N=\sum_{1\le i<j\le k}{\mathbb I}_{ij}$ denote
the number of disconnected pairs. Clearly, the event $N=0$ implies that all vertices from $\cup_{i\ge 1}C_i$ belong to the same connected component of $G$. Therefore, it suffices to show that $\PP(N= 0)=1-o(1)$. For this purpose we prove the
bound $\PP(N\ge 1|G')=o(1)$ uniformly in $G'$ satisfying (\ref{f-N}), see (\ref{2010-02-11A}) below.

In what follows we assume that (\ref{f-N}) holds. Let ${\hat f}(C_i)=\cup_{v\in C_i}S(v)$ denote the set of attributes occupied by vertices from $C_i$. Here
${\hat f}(C_i)\cap {\hat f}(C_j)=\emptyset$, for $i\not= j$. Note that if a red vertex finds neighbours
in $C_i$ and $C_j$ simultaneously
then it builds a path in $G$ that connects components $C_i$ and $C_j$. Clearly, only vertices with attribute sets
of size at least $2$ (i.e., vertices of types $2,3,\dots$) can build such a path.
The probability of building  such a path is minimized  by vertices of type $2$. This minimal probability is
\begin{equation}
\nonumber
p_{ij}
=2\frac{|{\hat f}(C_i)|\times|{\hat f}(C_j)|}{m(m-1)}.
\end{equation}
Note that (\ref{f-N}) combined with the inequality $|C_i|\ge \lceil n^{2/3}\rceil$ implies
$|{\hat f}(C_i)|\ge n^{2/3}(2M\ln n)^{-1}$. Hence,
\begin{equation}
\nonumber
p_{ij}
\ge \frac{1}{2M^2}\frac{n^{4/3}}{(m\ln n)^2}=:p_*.
\end{equation}
Let $r:=\lfloor n_2\varepsilon\rfloor+\dots+\lfloor n_M\varepsilon\rfloor$ denote the number of red vertices of types $2,3,\dots$. Observe that, for large $n$, (\ref{2010--A1}) implies $r\approx \varepsilon q'n$. Here $q'=q_2+\dots+q_M$.
In particular, we have
\begin{equation}\label{I-ij}
\PP(I_{ij}=1|G')\le (1-p_{ij})^r\le (1-p_*)^r\le e^{-p_*r}.
\end{equation}
Here $p_*r\ge c'n^{7/3}(\ln n)^{-2}$, and  the constant $c'$ depends on $\beta, M$, and  $q'$.
 Next, we apply  Markov's inequality to the conditional
probability
\begin{displaymath}
\PP(N\ge 1|G')\le \E(N|G')=\sum_{1\le i<j\le k}\PP(I_{ij}=1|G').
\end{displaymath}
Invoking (\ref{I-ij}) and the inequality $k\le (1-\varepsilon)n^{1/3}$ we obtain
\begin{equation}
\label{2010-02-11A}
\PP(N\ge 1|G')\le k^2e^{-p_*r}\le n^{2/3}e^{-c'n^{1/3}\ln^{-2}n}.
\end{equation}

\end{proof}

\begin{proof}[Proof of Lemma \ref{coupling}]
Throughout the proof we use the notation of Lemma \ref{Discrete}.

Fix $\omega(\cdot)$. Given $0<\varepsilon<1$,
let ${\cal Y}^{+\varepsilon}$  and ${\cal Y}^{-\varepsilon}$ be
multi-type Galton-Watson processes with type space $[M]$
where the number of children $Y^{+\varepsilon}_{ s t}$ ($Y^{-\varepsilon}_{ s t}$) of type $t$ of a particle
of type $s$ has binomial distribution
$Bi\bigl(\lfloor q_tn(1+\varepsilon)\rfloor, p_{st}(1+\varepsilon)\bigr)$ and
$Bi\bigl(\lfloor q_tn(1-\varepsilon)\rfloor, p_{st}(1-\varepsilon)\bigr)$ respectively.
Here $p_{st}:=(s-1)t(\beta n)^{-1}$.

Let ${\cal X}^{+\varepsilon}$ (and ${\cal X}^{-\varepsilon}$) be  multi-type Galton-Watson process with type space $[M]$
where the number of children  $X^{+\varepsilon}_{ s t}$ (and $X^{-\varepsilon}_{ s t}$) of type $t$ of a
particle of type $s$ has the Poisson distribution with  mean
 $\lambda_{st}(1+\varepsilon)$ (and $\lambda_{st}(1-\varepsilon)$). Here
$\lambda_{st}:=(s-1)tq_t\beta^{-1}$.

Given a multi-type  G-W process ${\cal Z}$ with type space $[M]$, by ${\cal Z}(t)$ we denote the process
starting at a particle of type $t$,
$|{\cal Z}(t)|$ denotes the total progeny of ${\cal Z}(t)$, $\rho({\cal Z},t):=\PP(|{\cal Z}(t)|=\infty)$ and
$\rho^{(k)}({\cal Z},t):=\PP(|{\cal Z}(t)|\ge k)$.

 It is known, see, e.g. inequality (1.23) in  \cite{Barbour},  that the total variation distance
 between the binomial
  distribution $Bi(r,p)$ and the Poisson distribution with the same mean is at most $p$.
  Therefore, by a coupling of the offspring numbers of binomial and Poisson branching processes we obtain
\begin{eqnarray}\label{2010-DD1}
&&
\rho^{(\omega(n))}({\cal Y}^{+\varepsilon},t)=\rho^{(\omega(n))}({\cal X}^{+\varepsilon'},t)+o(\omega(n)/n),
\\
\label{2010-DD2}
&&
\rho^{(\omega(n))}({\cal Y}^{-\varepsilon},t)=\rho^{(\omega(n))}({\cal X}^{-\varepsilon''},t)+o(\omega(n)/n).
\end{eqnarray}
Here $\varepsilon'=(1+\varepsilon)^2-1$ and $\varepsilon''=1-(1-\varepsilon)^2$. Letting  $n\to\infty$   we obtain,
\begin{equation}\label{2010-EE1}
\rho^{(\omega(n))}({\cal X}^{+\varepsilon'},t)\to\rho({\cal X}^{+\varepsilon'},t),
\qquad
\rho^{(\omega(n))}({\cal X}^{-\varepsilon''},t)\to\rho({\cal X}^{-\varepsilon''},t).
\end{equation}
Furthermore, letting $\varepsilon\downarrow 0$ we obtain
\begin{equation}
\label{2010-EE2}
\rho({\cal X}^{+\varepsilon'},t)\to\rho_{Q,\beta}(t),
\qquad
\rho({\cal X}^{-\varepsilon''},t)\to\rho_{Q,\beta}(t).
\end{equation}

{\it Proof of (\ref{2010-AA3}).}  We shall show that
\begin{equation}\label{2010-FF1}
\PP(v\in B^r)\le \rho_{Q,\beta}(s_v+1)+o(1).
\end{equation}
uniformly in $v\in V$. Collecting these bounds in the identity $\E|B^r|=\sum_{v\in V}\PP(v\in B^r)$ and using
(\ref{2010--A1})
we then obtain (\ref{2010-AA3}).
Therefore, it suffices to prove (\ref{2010-FF1}). In the proof we couple
  {\it regular} exploration starting at $v$ with the  process
${\cal Y}^{+\varepsilon}(s_v+1)$.
Let $Y^r_{it}$ denote the number of regular children of type $t$ discovered by $u_i\in L_v^r=\{v=u_1, u_2,\dots\}$.
Let $n_{it}$ denote the number of  uncoloured vertices of type $t$ at the moment, when
$u_i$ starts exploration of its neighbourhood.  Then  $Y^r_{it}$ has the binomial distribution
$Bi(n_{it}, p'_{it})$ with success probability
$p'_{it}=p_1(t,|{S'}^{r}(u_i)|, |W\setminus D_{i-1}|)$. Note that for large $n$ we have
\begin{equation}
\label{CCC-1}
n_{it}\le \lfloor q_tn(1+\varepsilon)\rfloor,
\qquad
p'_{it}\le |{S'}^{r}(u_i)|\, t(\beta n)^{-1}(1+\varepsilon).
\end{equation}
The first inequality follows from (\ref{2010--A1}). The second inequality follows from (\ref{pab})
combined with the inequalities
\begin{equation}
 m \ge |W\setminus D^r_{i-1}|=m-|D^r_{i-1}|> m-M\omega(n)=m-o(m).
\end{equation}
In addition, in view of (\ref{2010--A2}), we can replace $m$ by $\beta n$ in (\ref{CCC-1}).
(\ref{CCC-1}) shows that the parameters of
the binomial distribution of $Y^r_{it}$ are smaller than
the corresponding parameters of the offspring distribution of the branching process
${\cal Y}^{+\varepsilon}(s_v+1)$. Therefore,
particles of
 the branching process
produce at least as many children of each type as the vertices $u_i$, $i<\omega(n)$.
Note that  $v=u_1$ corresponds to a particle of type $|{S'}^r(v)|=s_v+1$ of the branching process while remaining
vertices
$u_i$, $i\ge 2$ correspond to particles of types $s_{u_i}=|S(u_i)|$ respectively.
Hence, we have
\begin{equation}\label{2010-FF3}
\PP(v\in B^r)\le \PP\bigl(|{\cal Y}^{+\varepsilon}(s_v+1)|\ge \omega(n)\bigr).
\end{equation}
(\ref{2010-FF3}) in combination with
(\ref{2010-DD1}), (\ref{2010-EE1}) and (\ref{2010-EE2}) implies (\ref{2010-FF1}).


 {\it Proof of (\ref{2010-AA1}).}
Given $v\in V$, we start simple exploration at $v$. Let $K_{t}$ ($I_{t}$) denote the  number of complex
(irregular) children of type $t$
discovered  by the exploration until the list $L_v^s(\omega(n))$
was completed. We put a label on $v$ whenever $\max_t\{K_{t},I_{t}\}\ge \omega(n)$.

Let $A$ denote the set of labeled vertices and $p_v':=\PP(v\in B^s|\,v\notin A)$ be
 the probability that the simple exploration of  unlabeled vertex $v$ discovers
at least $\omega(n)$ vertices. We show below that
\begin{eqnarray}\label{2010-GG001}
&&
\PP(v\in A)
=O(n^{-1}),
\\
\label{2010-HH001}
&&
p_v'
=\rho_{Q,\beta}(s_v+1)-o(1).
\end{eqnarray}

It follows from (\ref{2010-GG001}), (\ref{2010-HH001})  that
\begin{equation}\label{bzzz}
\PP(v\in B^s)=p_v'+O(n^{-1})=\rho_{Q,\beta}(s_v+1)+o(1).
\end{equation}
Invoking the latter identity in the expression $\E|B^s|=\sum_{v\in V}\PP(v\in B^s)$ we obtain
(\ref{2010-AA1}).

\bigskip

{\it Proof of (\ref{2010-HH001}).}
Given $\varepsilon>0$ we show that for large $n$
\begin{equation}\label{2010-hg}
\PP\bigl(|{\cal Y}^{+\varepsilon}(s_v+1)|\ge \omega(n)\bigr)
\ge
p_v'
\ge
\PP\bigl(|{\cal Y}^{-\varepsilon}(s_v+1)|\ge \omega(n)\bigr).
\end{equation}
These inequalities in combination with (\ref{2010-DD1}-\ref{2010-EE2}) imply (\ref{2010-HH001}).

 In order to generate events of probability $p'_v$ we use rejection  sampling.
In the course of exploration we keep track of the number of coloured vertices and interrupt the exploration at the moment when this number exceeds $3\omega(n)$. Exploration is rejected if it is interrupted before the list $L_v^s(\omega(n))$ is completed. Otherwise it is accepted. Clearly, $p'_v$ is the probability
that the list $L_v^s(\omega(n))$ of an accepted exploration has collected all $\omega(n)$ elements.

In the proof of
(\ref{2010-hg}) we couple
the simple exploration process with branching processes ${\cal Y}^{-\varepsilon}(s_v+1)$ and
${\cal Y}^{+\varepsilon}(s_v+1)$ so that the number of simple children of type $t$ of the vertex $v$ is at least (most)
 as large
as the number of particles of type $t$ in the first generation of ${\cal Y}^{-\varepsilon}(s_v+1)$
(${\cal Y}^{+\varepsilon}(s_v+1)$),  $t\in [M]$.
In the further steps of exploration the number $Y_t(u)$ of simple children of type $t$ discovered by a particle
$u\in L_v^s(\omega(n))\setminus \{v\}$ is at least (most) as large as the number of children of type $t$ produced by the coresponding particle of type $s_u$ of the process ${\cal Y}^{-\varepsilon}$ (${\cal Y}^{+\varepsilon}$).

To make sure that such a coupling is possible we fix $u=u_i\in L_v^s(\omega(n))$ and count its simple children.
Recall that $u_i$ selects simple children from the current set of uncoloured vertices. These are checked one after another in increasing order, and
each newly discovered simple child is added to the list $L_v^s$ before the next uncoloured vertex is checked.
At the moment when a vertex  $g$ is checked,
its probability
to be a simple child of $u$ is $p_i(g)=p(|S(g)|, |{S'}^s(u)|, |H_i(g)|, |W \setminus D_{i-1}|)$.
It is a conditional probability given $\{S(u'),\, u'\in L_v^s\}$. Here $L_v^s$ is the set of
vertices that have been added to the list before $g$ was checked.  Note that,
as far as the probability of the event $\{v\in B^s\}\equiv \{L_v^s(\omega(n))=\omega(n)\}$ is considered,
 we may safely assume that
$|D_{i-1}|, |H_i(g)|\le M(\omega(n)-1)$.
 It follows from these inequalities and (\ref{pabh}) that for large $n$ we have
\begin{equation} \label{01-26-1}
\frac{|{S'}^s(u)|s_g}{m}(1-\varepsilon)\le
 p_i(g)
\le
\frac{|{S'}^s(u)|s_g}{m}(1+\varepsilon).
\end{equation}
In addition, in view of (\ref{2010--A2}), we can replace $m$ by $\beta n$ in the denominator.
Let $n_{it}^*$ denote the number of uncoloured vertices of type $t$ at the moment when $u=u_i$
starts search of its simple children. Until the exploration is not interrupted we have
$n_{it}^*\ge n_t-3\omega(n)$. For large $n$ this inequality
implies $n_{it}^*\ge (1-\varepsilon/2)n_t$. Invoking (\ref{2010--A1}) we obtain
\begin{equation} \label{01-26-2}
q_tn(1-\varepsilon)\le n_{it}^*\le q_tn(1+\varepsilon)
\qquad
t\in [M].
\end{equation}
It follows from (\ref{01-26-1}, \ref{01-26-2}) that we can couple $Y_t(u)$ with binomial random variables
\begin{displaymath}
Y^{\pm}_t(u)\sim Bi\Bigl(\lfloor q_tn(1\pm\varepsilon)\rfloor,\frac{|{S'}^s(u)|s_g}{\beta n}(1\pm\varepsilon)\Bigr),
\end{displaymath}
 so that almost surely we have
$Y^{-}_t(u)\le Y_t(u)\le Y^{+}_t(u)$. These inequalities imply (\ref{2010-hg}).

{\it Proof of (\ref{2010-GG001})}.
We write $ \PP(v\in A)\le \sum_{t\in [M]}(\PP(K_t\ge\omega(n))+\PP(I_t\ge \omega(n))$
  and show that
 \begin{equation}\label{2010-GG-001A}
 \PP(K_t\ge\omega(n))=o(n^{-1}),
 \qquad
 \PP(I_t\ge\omega(n))=o(n^{-1}).
 \end{equation}
 We prove the first bound only. The proof of the second bound is much the same.
Given $i\le \omega(n)$,
the number of complex children of type $t$ discovered by
 $u_i\in L_v^s$
is the sum of at most $n_t$ independent Bernoulli random variables each with success probability at most
\begin{displaymath}
p^* =p_1\bigl(M,M,M\omega(n), m-M\omega(n)\bigr)\le cM^4m^{-2},
\end{displaymath}
see (\ref{pabh+1}).
Therefore,  $K_{t}$ is at most  sum of $n_t\omega(n)$ independent Bernoulli random variables
  with success probability $p^*$. In particular, we have
  \begin{equation}\label{2010-GG-001B}
\PP(K_t\ge\omega(n))\le\PP(\xi\ge\omega(n)),
\end{equation}
where $\xi\sim Bi(n_t\omega(n), p^*)$.
 By Chebychev's inequality
\begin{equation}\label{2010-GG-001C}
\PP(\xi\ge\omega(n))\le (\omega(n)-\E \xi)^{-2}\Var \xi =O(n^{-1}).
\end{equation}
In the last step we invoke the simple bounds
\begin{displaymath}
\Var \xi\le \E\xi=n_t\omega(n)p^*_t =O(\omega^2(n)n^{-1})=o(\omega(n)).
\end{displaymath}
(\ref{2010-GG-001B}) and (\ref{2010-GG-001C}) imply the first bound of (\ref{2010-GG-001A}).

{\it Proof of (\ref{2010-A2})}.  It suffices to establish (\ref{2010-A2}) for one particular function $\omega$,
because for any
other ${\tilde B}^s$
defined by another such function ${\tilde \omega}$,  we have
\begin{equation}\label{BB1}
|B^s|-|{\tilde B}^s|=o_P(n).
\end{equation}
To see this
write
$|B^s|-|{\tilde B}^s|\le |B^s\cup {\tilde B}^s|-|B^s\cap {\tilde B}^s|$ and
observe that $B^s\cup {\tilde B}^s$ and $B^s\cap {\tilde B}^s$
represent sets of bs-vertices defined by the functions $\omega_1=\min\{\omega,{\tilde \omega}\}$ and
$\omega_2=\max\{\omega,{\tilde \omega}\}$ respectively.
An application of (\ref{2010-AA1})  to $\omega_1$ and $\omega_2$ yields the bound
$\E(|B^s|-|{\tilde B}^s|)=o(n)$. This bound implies
(\ref{BB1}).

We show (\ref{2010-A2}) for $\omega(n)=\lfloor \ln n \rfloor$. For this purpose we prove the bound for the variance
\begin{equation}\label{Var}
\E|B^s|^2-(\E|B^s|)^2=o(n^2),
\end{equation}
which tells us that $|B^s|-\E|B^s|=o_P(n)$.
In particular, (\ref{Var})
combined with (\ref{2010-AA1})
 shows (\ref{2010-A2}).

In the proof of (\ref{Var}) we use the observation that the first
$\omega(n)$ steps of
any two explorations starting at distinct vertices are almost independent.
More precisely, we show below that uniformly in $\{u,v\}\subset V$
\begin{equation}\label{Var-1}
\PP(u,v\in B^s)=\rho_{Q,\beta}(s_u+1)\rho_{Q,\beta}(s_v+1)+o(1).
\end{equation}
It follows from (\ref{Var-1}) that
\begin{eqnarray}
\label{Var-2}
2\sum_{\{u,v\}\subset V}\PP(u,v\in B^s)
&
=
&
\sum_{u,v\in V}\rho_{Q,\beta}(s_u+1)\rho_{Q,\beta}(s_v+1)+o(n^2).
\\
\nonumber
&
=
&
n^2{\tilde\rho}^2+o(n^2).
\end{eqnarray}
In the last step we use (\ref{2010--A1}).
Observe, that the left-hand sum of (\ref{Var-2}) is the expected value of
$2\sum_{\{u,v\}\subset V}{\mathbb I}_{\{u, v\in B^s\}}=|B^s|^2-|B^s|$. Therefore, from (\ref{Var-2})
we obtain
\begin{displaymath}
\E|B^s|^2=n^2{\tilde \rho}^2+\E|B^s|+o(n^2).
\end{displaymath}
 This identity combined with
 (\ref{2010-AA1})
implies (\ref{Var}).

\medskip

Let us prove (\ref{Var-1}). We  first explore  $u$ and then $v$. In each case we stop simple exploration after the number
of vertices in the corresponding list
reaches $\omega(n)$.
 Note that with a high probability these
two explorations do not meet.
Indeed, let $T_u$ ($T_v$) denote the set of vertices coloured by the first (second)  exploration and
let ${\cal H}$ denote the event that the second exploration does not encounter any vertex from $T_u$,  i.e.,
${\cal H}=\{D_u\cap S(v')=\emptyset$, for each $v'\in T_v\}$. Here we denote $D_u=\cup_{u'\in T_u}S(u')$ and
$D_v=\cup_{v'\in T_v}S(v')$.
Now assume that $u,v$ are  unlabeled vertices, i.e.,  $u,v\notin A$. Then
\begin{displaymath}
|T_u|, |T_v|\le (2M+1)\omega(n)=:{\hat T},
\end{displaymath}
and
$|D_u|, |D_v|< M {\hat T}\le M(2M+1)\omega(n)=:{\hat D}$.
In this case, for each $v'\in T_v$, the probability that  $S(v')$ does not hit $D_u$ is at least
$\bigl(\frac{m-2{\hat D}}{m}\bigr)^{M}$. Here we use the fact that
 $S(v')$ has at most $|S(v')|\le M$ elements (trials) to hit the set
$D_u$ which occupies $|D_u|\le {\hat D}$ attributes  among those (at least $m-{\hat D}$) which have
not been used  by
the current collection of  vertices of evolving list $L_v^s$. Since there are at most ${\hat T}$ vertices in $T_v$, we obtain
\begin{displaymath}
\PP({\cal H}|u,v\notin A)
\ge
\bigl(\frac{m-2{\hat D}}{m}\bigr)^{M{\hat T}}
=1-O\bigl(\omega^2(n)n^{-1}\bigr).
\end{displaymath}
For arbitrary $u,v$ we obtain from (\ref{2010-GG001})
\begin{equation}\label{2010-HS}
\PP({\cal H})\ge \PP({\cal H}\cap\{u,v\notin  A\})
=
\PP({\cal H}|u,v\notin A)\PP(u,v\notin A)=1-o(1).
\end{equation}
Now assume that $\rho_{Q,\beta}(s_u+1)>0$  (otherwise (\ref{Var-1})  trivially follows from (\ref{bzzz})) and write
\begin{equation}\label{puga-3}
\PP(u,v\in B^s)
=
\PP(v\in B^s|u\in B^s)\PP(u\in B^s).
\end{equation}
  We can
replace
$\PP(v\in B^s|u\in B^s)$
by $p_{v,u}:=\PP\bigl(v\in B^s\bigr|\{u\in B^s\}\cap\{u,v\in A\}\cap {\cal H}\bigr)$
and $\PP(u\in B^s)$ by  $\rho_{Q,\beta}(s_u+1)$. It follows from (\ref{2010-GG001}), (\ref{2010-HS}) and
(\ref{bzzz}) that the error due to such replacement is of order $o(1)$. From (\ref{puga-3}) we obtain
\begin{equation}\label{puga0102}
\PP(u,v\in B^s)=p_{v,u}\rho_{Q,\beta}(s_u+1)+o(1).
\end{equation}
Finally, (\ref{Var-1}) follows from (\ref{puga0102}) and the identity
$p_{v,u}=\rho_{Q,\beta}(s_v+1)+o(1)$,
which is shown in much the same way as (\ref{2010-HH001}) above.
%
\end{proof}

\begin{proof}[Proof of Lemma \ref{inequalities-1}]
 Let $(x_1,\dots, x_k)$ be a random permutation of elements of the set  $W'$. For $A=\{x_1,\dots, x_a\}$ we have,
by symmetry,
\begin{eqnarray}
\label{pab+1}
p(a,b,k)
&&
\le \sum_{1\le i\le a}\PP(x_i\in B)=a\PP(x_1\in B),
\\
\label{p1ab+1}
p_1(a,b,k)
&&
=\sum_{1\le i\le a}\PP(A\cap B=x_i)=a\PP(A\cap B=x_1),
\\
\label{p2abm+1}
p_2(a,b,k)
&&
\le \sum_{1\le i<j\le a}\PP(x_i,x_j\in B)=2^{-1}a(a-1)\PP(x_1,x_2\in B),
\\
\label{pabh+11}
p(a,b,h,k)
&&
=\sum_{1\le i\le a}\PP(A\cap B=x_i)\PP(H\cap A=\emptyset|A\cap B=x_i)
\\
\nonumber
&&
=p_1(a,b,k)\bigl(1-p(a-1,h,k-b)\bigr).
\end{eqnarray}
The right-hand side inequality of (\ref{pab}) follows from (\ref{pab+1}) and the identity $\PP(x_1\in B)=b/k$.
The left-hand side inequality follows from
 (\ref{p1ab+1}) combined with the identity
$\PP(A\cap B=x_1)=\frac{b(k-b)_{a-1}}{(k)_a}$ and inequalities
\begin{displaymath}
1\ge \frac{(k-b)_{a-1}}{(k-1)_{a-1}}\ge \bigl(\frac{k-a-b}{k-a}\bigr)^{a-1}\ge 1-\frac{ab}{k-a}.
\end{displaymath}
(\ref{p2ab}) follows from  (\ref{p2abm+1}) and the identity
$\PP(x_1,x_2\in B)=\frac{(b)_2}{(k)_2}$.
(\ref{pabh}) follows from (\ref{pabh+11}) combined with (\ref{pab}).
(\ref{pabh+1}) follows from the inequality $p_1(a,b,h,k)=p_1(a,b,k)p(a-1,h,k-b)$, which is shown in the same way as
(\ref{pabh+11}).
 \end{proof}

\end{document}